\documentclass[11pt,a4paper]{amsart}
\usepackage[T1]{fontenc}
\usepackage{lmodern,amssymb}
\usepackage{graphicx,float}
\usepackage[pass]{geometry}
\usepackage[hidelinks]{hyperref}

\hypersetup{
  pdftitle={The Koebe Radius for Univalent Cubic Polynomials},
  pdfauthor={Alexey Solyanik}
}

\newtheorem{theorem}{Theorem}

\title[The Koebe radius for univalent cubics]
      {The Koebe Radius for Univalent Cubic Polynomials}
\author{Alexey Solyanik}
\address{Odesa Polytechnic National University, Odesa, Ukraine}
\email{solyanik@op.edu.ua}
\date{November 14, 2026}

\begin{document}
\raggedbottom

\begin{abstract}
Let \(p(z)=z+a_2z^2+a_3z^3\) be a polynomial univalent in the unit disk
\(\Delta\). We determine the largest disk centered at the origin that
is contained in the image \(p(\Delta)\) of every such polynomial. Its
squared radius is the unique real root of
\[
 27q^3-54q^2+63q-8=0.
\]
We also find the extremal polynomial \(p_*\). Quite surprisingly, the
image \(p_*(\Delta)\) of the extremizer has no symmetry.
\end{abstract}

\maketitle

Let \(\mathcal S\) denote the class of functions that are analytic and
univalent (one-to-one) in the unit disk \(\Delta\) of the complex plane,
normalized by the conditions \(f(0)=0\), \(f'(0)=1\). Thus
\begin{equation}\label{eq:normalized-functions}
 f(z)=z+a_2z^2+a_3z^3+\cdots .
\end{equation}

In his 1907 paper on the uniformization of analytic curves, Koebe
proved that there is a number \(r>0\), independent of \(f\), such that
\[
 \Delta_r=\{z\in\mathbb C:|z|<r\}
\]
is contained in \(f(\Delta)\) for every \(f\in\mathcal S\)
\cite{Koebe}. This estimate was needed to prevent normalized univalent
maps from degenerating in a limiting construction for the general
uniformization problem \cite{Koebe}. Koebe conjectured that the exact
value is \(r=1/4\) and that an extremal function is
\begin{equation}\label{eq:koebe-function}
 k(z)=\frac{z}{(1-z)^2}=z+2z^2+3z^3+\cdots ,
\end{equation}
which maps \(\Delta\) conformally onto the slit plane
\(\mathbb C\setminus(-\infty,-1/4]\).

Bieberbach proved this result in 1916 by deriving it from the sharp
coefficient estimate \(|a_2|\le2\) \cite{Bieberbach} (see also
\cite[Chapter~2]{Duren}). In the same paper, Bieberbach formulated his
celebrated conjecture \(|a_m|\le m\), \(m=2,3,\ldots\), for every
function in \(\mathcal S\), with the Koebe function as an extremal
function.

This conjecture shaped much of the subsequent development of univalent
function theory in the twentieth century and was finally proved by
Louis de Branges. His proof was published in 1985 \cite{deBranges}.

To the best of our knowledge, univalent polynomials were first studied
systematically by Alexander in his fascinating 1915 paper
\cite{Alexander}. He related univalence to the locations of the zeros
of a polynomial and its derivative and established several sufficient
conditions for univalence.

Much of the subsequent work on normalized
univalent polynomials concentrated on the coefficient problem, in the
spirit of the Bieberbach conjecture. In particular, in 1969 Suffridge
\cite{Suffridge} constructed a family of polynomials from \(\mathcal S\) whose
coefficients tend to those of the Koebe function
\eqref{eq:koebe-function}.

The general problem of describing the coefficient regions of univalent
polynomials appears to be extremely difficult and seems unlikely to
admit a tractable description. For polynomials of low degree, however,
such a description is possible. In particular, Brannan in 1967
\cite{Brannan}, and independently Cowling and Royster in 1968
\cite{CowlingRoyster}, described all complex pairs \((a_2,a_3)\) for
which
\[
 p(z)=z+a_2z^2+a_3z^3
\]
is univalent; see also the later account of Rahman and Waniurski
\cite[Introduction]{RahmanWaniurski}.

A normalized rotation of \(p\) is the polynomial
\(p_\theta(z)=e^{-i\theta}p(e^{i\theta}z)\), where \(\theta\) is real.
We have \(p_\theta'(0)=1\) and \(p_\theta(\Delta)=e^{-i\theta}p(\Delta)\),
so the largest origin-centered disk in the image has the same radius. The cubic
coefficient becomes \(e^{2i\theta}a_3\). Choosing \(\theta\) so that
this coefficient equals \(b=|a_3|\), we may write
\begin{equation}\label{eq:normalized-cubic}
 p(z)=z+az^2+bz^3,\qquad a=u+iv,
\end{equation}
where \(u\) and \(v\) are real.
Then Brannan's criterion \cite{Brannan} says that \(p\) is univalent
in \(\Delta\) if and only if
\begin{equation}\label{eq:intro-b-range}
 0\leq b\leq\frac13,
\end{equation}
and, for every \(s\in[0,4]\) with \(s=4\) omitted when \(b=1/3\),
\begin{equation}\label{eq:intro-ellipses}
 \frac{s u^2}{(1+b(s-1))^2}
 +\frac{s v^2}{(1-b(s-1))^2}\leq1.
\end{equation}

\begin{figure}[H]
\centering
\makebox[\textwidth][c]{%
  \hspace*{-0.7cm}%
  \includegraphics[width=1.08\textwidth]{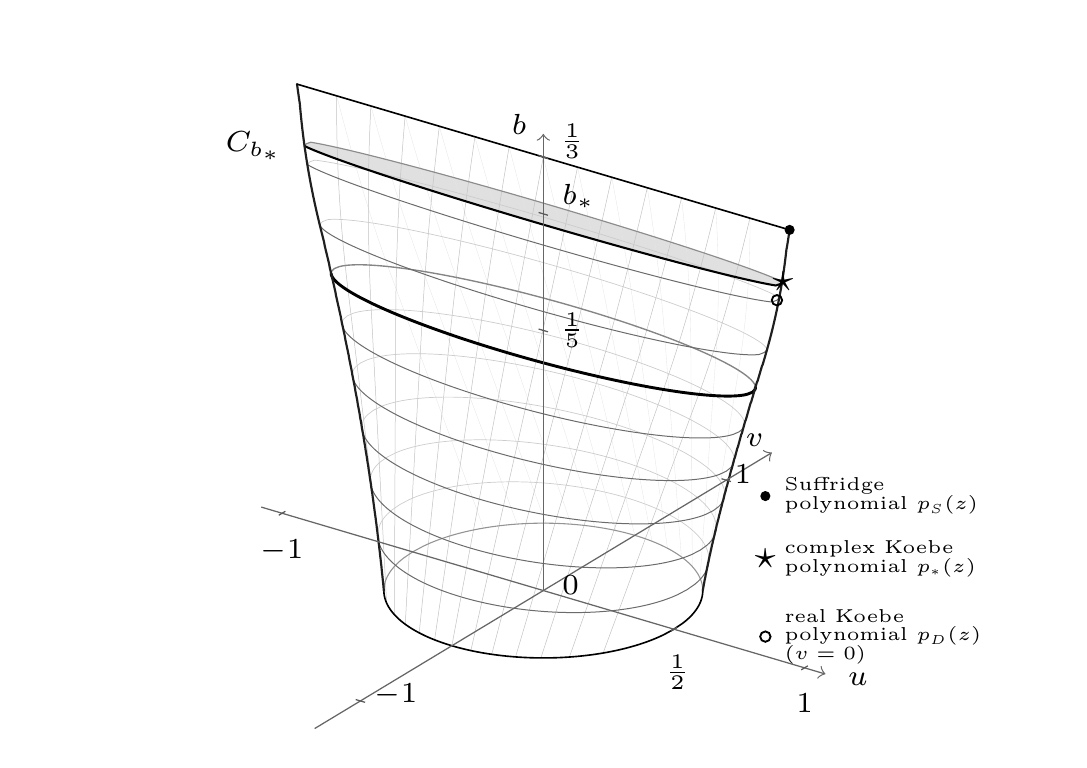}%
  \hspace*{0.7cm}%
}
\caption{The coefficient body of normalized univalent cubic
polynomials.}
\label{fig:coefficient-body}
\end{figure}

\begin{figure}[H]
\centering
\includegraphics[width=.82\textwidth]{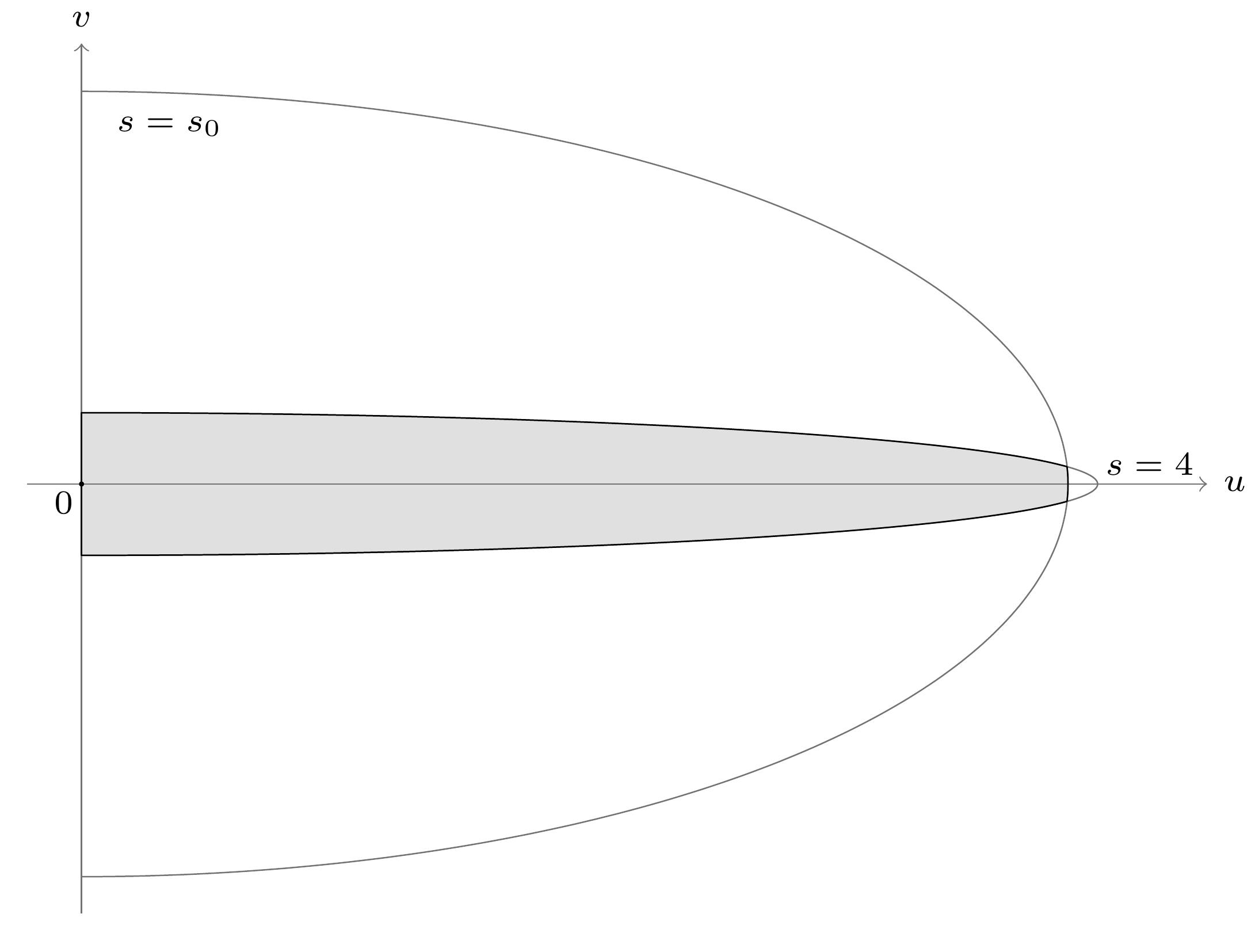}
\caption{One half of the section \(C_b\) of the coefficient body at
\(b=b_*\approx0.2896\).}
\label{fig:coefficient-section}
\end{figure}

For each \(b\) satisfying \eqref{eq:intro-b-range}, let \(C_b\) be
the region of admissible coefficients \(a=u+iv\). It is the intersection
of the ellipses in \eqref{eq:intro-ellipses} for \(0<s\le4\), with
\(s=4\) omitted when \(b=1/3\). At \(s=0\) the inequality holds
identically. A direct calculation shows
that, for \(0\leq b\leq1/5\), this intersection is itself an ellipse;
when \(b=0\), it is a circle. For \(1/5<b<1/3\), no single ellipse
determines the intersection: its boundary consists of an arc of the
ellipse corresponding to \(s=4\) and the envelope of a nontrivial
subfamily of \eqref{eq:intro-ellipses}, as in
Figure~\ref{fig:coefficient-section}. When \(b=1/3\), letting
\(s\uparrow4\) forces \(v=0\); the remaining inequalities give the
segment \(|u|\le2\sqrt2/3\).

Brannan's criterion immediately gives
\[
 |a_2|\leq\frac{2\sqrt2}{3},\qquad |a_3|\leq\frac13,
\]
and the cubic Suffridge polynomial
\[
 p_S(z)=z+\frac{2\sqrt2}{3}z^2+\frac13z^3
\]
maximizes both coefficient moduli
\cite{Brannan,CowlingRoyster,Suffridge}.
Thus the coefficient problem for cubic polynomials was completely
solved.

To the best of our knowledge, Dimitrov first posed a polynomial analogue
of Koebe's one-quarter theorem in 2002 \cite{Dimitrov}. Namely, for each
\(n\), determine the largest disk centered at the origin that is
contained in \(p(\Delta)\) for every polynomial \(p\in\mathcal S\) of
degree \(n\). He suggested that the Suffridge polynomials might be
extremal. We call the extremizers in Dimitrov's problem
\emph{Koebe polynomials}.

For a Suffridge polynomial \(p_S\), let \(r_S\) denote the radius of
the largest origin-centered disk contained in \(p_S(\Delta)\). For the
cubic Suffridge polynomial above, a direct calculation gives
\[
 r_S=\frac{2}{3\sqrt3}.
\]

\begingroup
\setlength{\intextsep}{4pt}
\begin{figure}[H]
\centering
\makebox[\textwidth][c]{%
  \includegraphics[height=190bp]{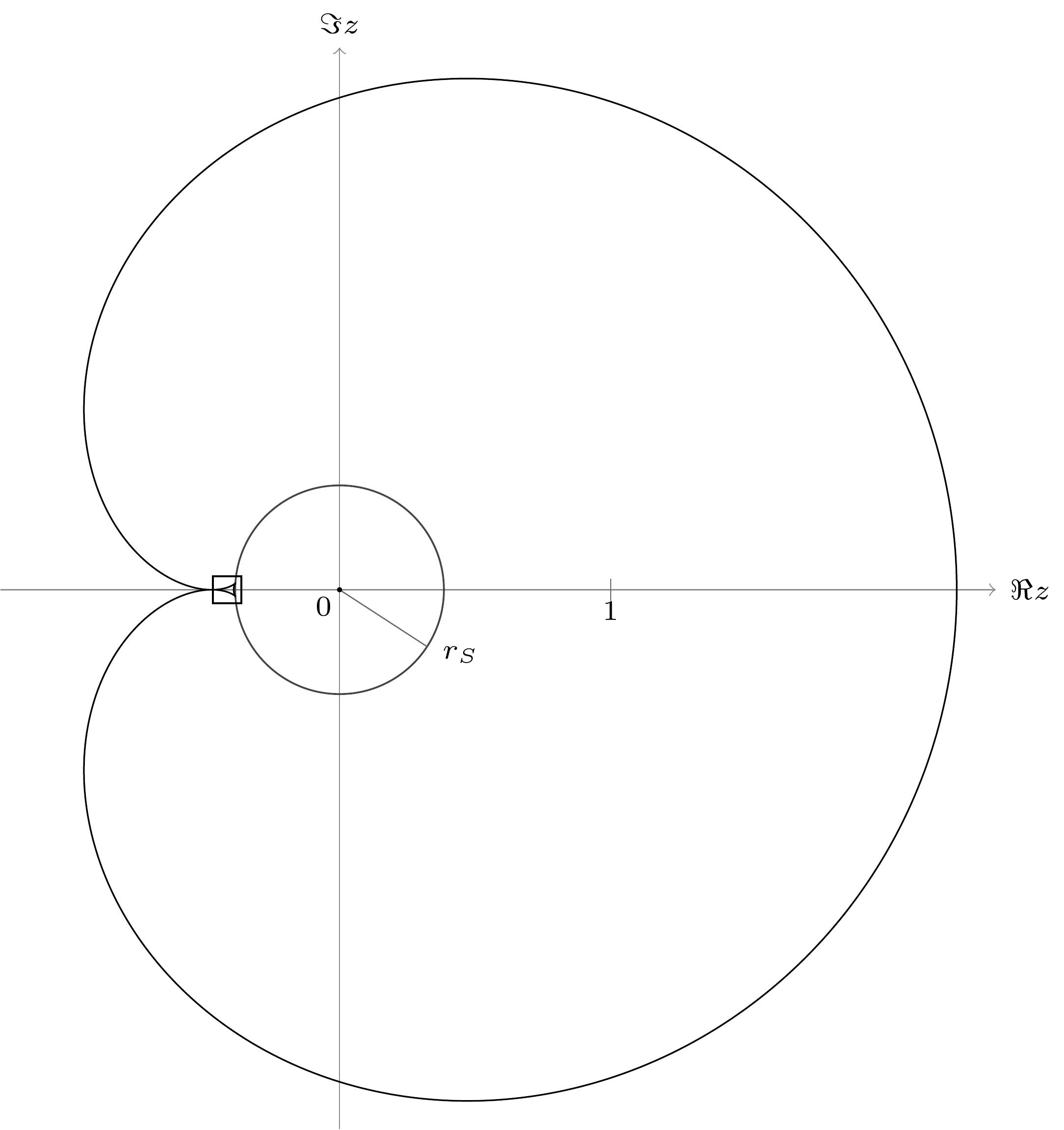}%
  \hspace{10bp}%
  \includegraphics[height=190bp]{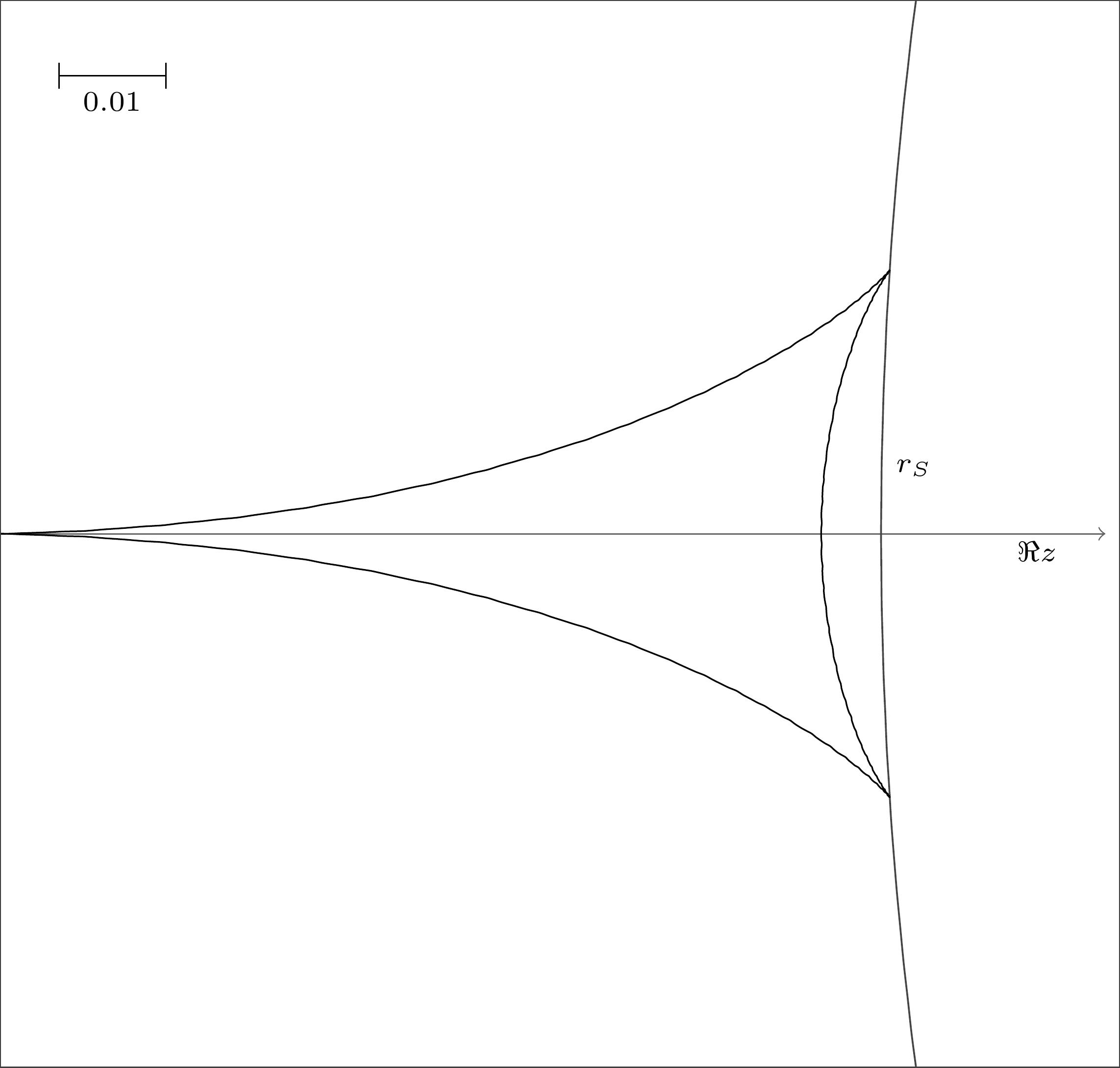}%
}
\caption{The image \(p_S(\Delta)\) of the cubic Suffridge polynomial
and a magnification of the boxed region.}
\label{fig:suffridge-image}
\end{figure}

In 2019, Dmitrishin, Dyakonov and Stokolos proved
\cite[Theorem~5.1]{DDS} that the Koebe radius for
\emph{cubic polynomials from \(\mathcal S\) with real coefficients} is
\[
 r_D=\frac{3-\sqrt5}{2},
\]
and that a Koebe polynomial in this subclass is
\[
 p_D(z)=z+\frac{2}{\sqrt5}z^2
   +\frac12\left(1-\frac1{\sqrt5}\right)z^3.
\]

\begin{figure}[H]
\centering
\makebox[\textwidth][c]{%
  \includegraphics[height=190bp]{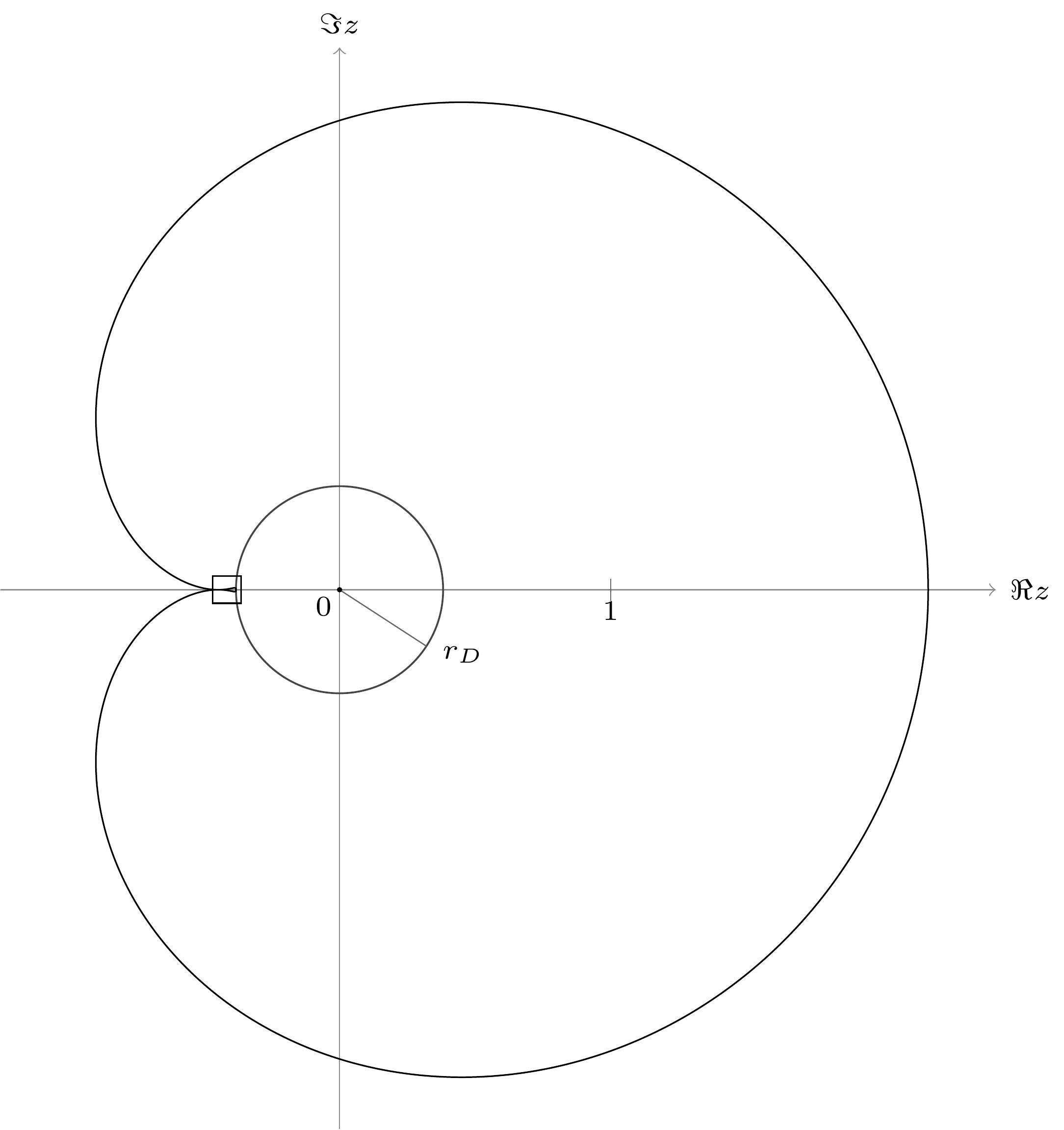}%
  \hspace{10bp}%
  \includegraphics[height=190bp]{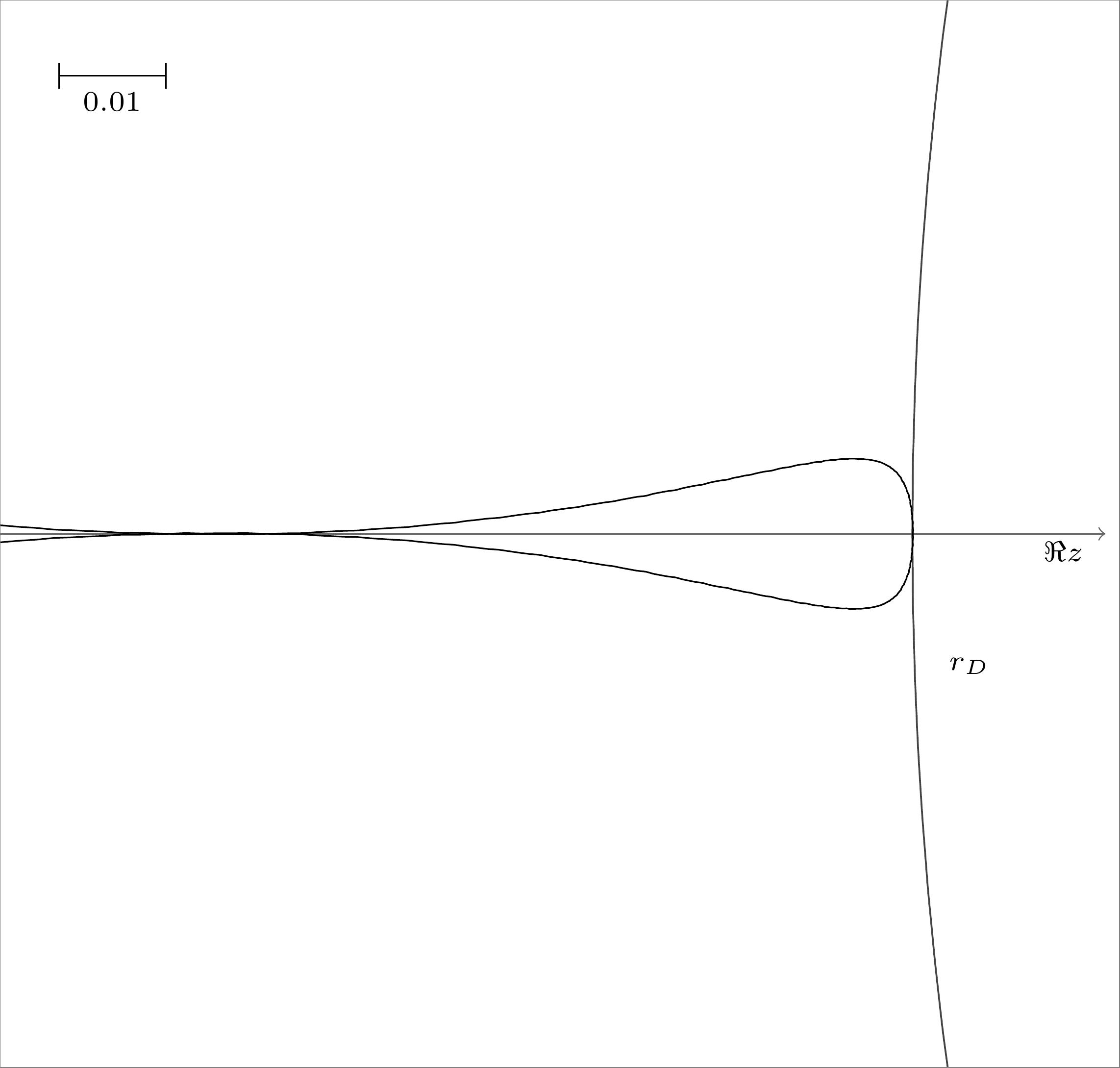}%
}
\caption{The image \(p_D(\Delta)\) of the real Koebe polynomial and a
magnification of the boxed region.}
\label{fig:real-cubic-image}
\end{figure}

Since \(r_D<r_S\), this result shows that the cubic Suffridge polynomial
is not a Koebe polynomial even within the subclass of polynomials with
real coefficients. The Koebe problem for the full class of univalent
cubic polynomials with complex coefficients remained open.

In this paper we determine the Koebe radius \(r_*\) for cubic
polynomials from \(\mathcal S\) and construct a cubic polynomial \(p_*\)
attaining it.

\begin{figure}[H]
\centering
\makebox[\textwidth][c]{%
  \includegraphics[height=190bp]{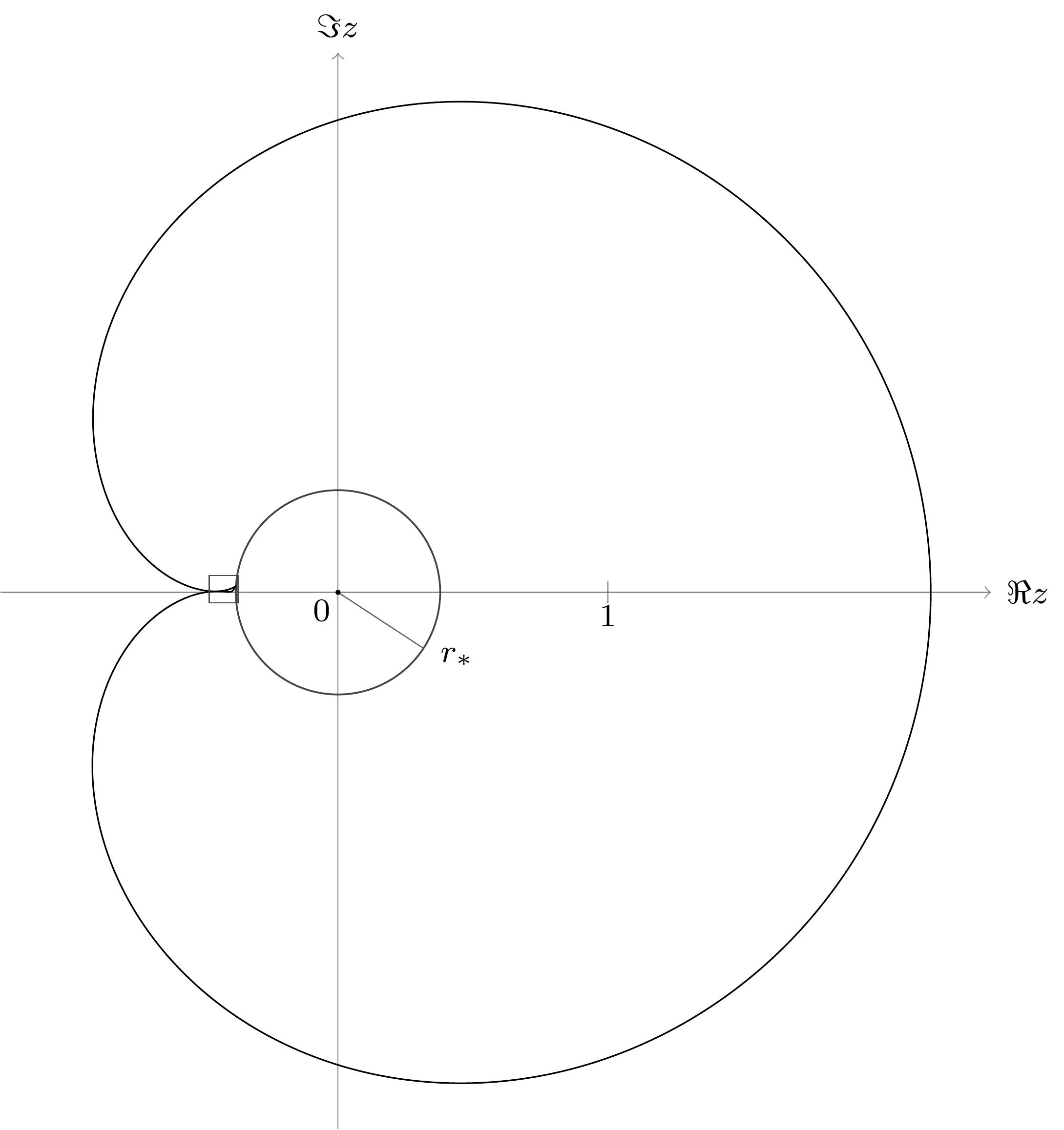}%
  \hspace{10bp}%
  \includegraphics[height=190bp]{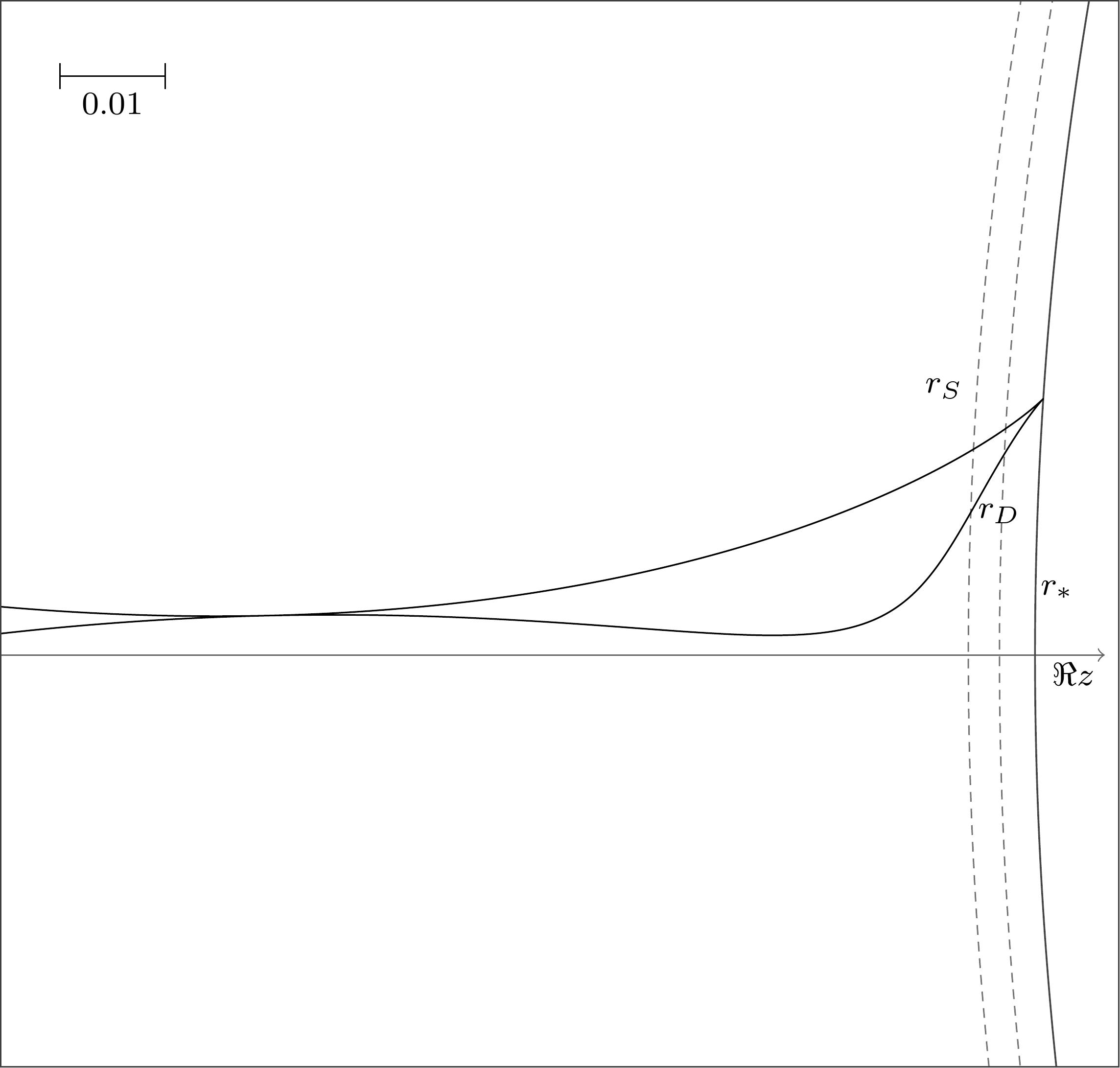}%
}
\caption{The image \(p_*(\Delta)\) and a magnification of the boxed
region.}
\label{fig:cubicimage}
\end{figure}
\endgroup

The exact value of \(r_*\) is given in Theorem~\ref{thm:main}. In
particular, \(r_S>r_D>r_*\), where \(r_*\) is the Koebe radius for
univalent cubic polynomials. One critical point of \(p_*\) lies on
\(\partial\Delta\), while the other lies outside \(\Delta\). The
image of the boundary critical point is a cusp of
\(p_*(\partial\Delta)\) and lies on \(\partial\Delta_{r_*}\).

Before stating and proving the theorem, we describe the main ideas
that lead to the solution. For each univalent cubic \(p\), the largest disk centered at the
origin and contained in \(p(\Delta)\) has radius
\[
 r(p)=\min_{|z|=1}|p(z)|.
\]
For \(p\) in the form \eqref{eq:normalized-cubic}, we have
\(p(z)/z^2=a+\overline z+bz\) on the unit circle. The curve traced by
\(\overline z+bz\) is centrally symmetric, so
\begin{equation}\label{eq:overview-curve-distance}
 r(p)=\min_{|z|=1}\bigl|a-(\overline z+bz)\bigr|.
\end{equation}
Since \(\overline z+bz\), as \(z\) ranges over the unit circle,
describes the boundary of the ellipse
\begin{equation}\label{eq:outer-ellipse}
 E_b=\left\{\xi+i\eta\in\mathbb C:
 \frac{\xi^2}{(1+b)^2}+\frac{\eta^2}{(1-b)^2}\le1\right\},
\end{equation}
we obtain
\begin{equation}\label{eq:overview-distance-to-ellipse}
 r(p)=\operatorname{dist}(a,\partial E_b).
\end{equation}
For fixed \(b\), we therefore seek the point of \(C_b\) nearest to
\(\partial E_b\), and then minimize the resulting distance over \(b\).

Thus finding the sharp Koebe radius and the extremal polynomial reduces
to minimizing the distance between the region \(C_b\) of admissible second
coefficients and the boundary of the ellipse \(E_b\) determined by the
third coefficient \(b\).

For \(0\le b\le1/4\), the ellipse \(K_b\) corresponding to \(s=3\)
in Brannan's family contains \(C_b\). Comparing \(K_b\) with
\(\partial E_b\) gives the uniform estimate
\begin{equation}\label{eq:overview-first-range}
 r(p)\ge\sqrt{\frac7{48}}.
\end{equation}
This lower bound already exceeds the Koebe radius attained in the
second range.

For \(1/4\le b<1/3\), we obtain the exact distance. Both estimates
apply at \(b=1/4\). Two ellipses in Brannan's family, corresponding
to \(s=s_0(b)\) and \(s=4\), meet at a corner \(a_b\) of \(C_b\).
A positive linear combination of their defining quadratic forms gives
an ellipse \(L_b\) that contains \(C_b\) and passes through \(a_b\).
The shortest segment from \(L_b\) to \(\partial E_b\) can be taken
to start at \(a_b\), so
\[
 \operatorname{dist}(L_b,\partial E_b)
 =\operatorname{dist}(C_b,\partial E_b)
 =\operatorname{dist}(a_b,\partial E_b).
\]
The case \(b=1/3\) follows by a limit. Minimizing the distance over
\(b\) gives the extremal parameter \(b_*\) and the radius \(r_*\).

\begin{theorem}
\label{thm:main}
Let
\[
 p(z)=z+a_2z^2+a_3z^3,\qquad a_2,a_3\in\mathbb C,
\]
be univalent in \(\Delta\). Then
\[
 \Delta_{r_*}\subset p(\Delta),
\]
where \(r_*^2\) is the unique real root of
\[
 27q^3-54q^2+63q-8=0.
\]
Namely,
\begin{equation}\label{eq:r-star}
 r_*=
 \sqrt{\frac{2+\sqrt[3]{6\sqrt3-9}-\sqrt[3]{6\sqrt3+9}}{3}}.
\end{equation}
The radius is sharp: there exists an extremizer \(p_*\) such that
\(\Delta_\rho\not\subset p_*(\Delta)\) for every \(\rho>r_*\). Up to
normalized rotations and complex conjugation, this extremizer is
unique and is given by
\[
 p_*(z)=z+(u_*+iv_*)z^2+b_*z^3,
\]
where
\[
 \begin{aligned}
 b_*&=\sqrt{r_*^4-\frac{r_*^2}{3}+\frac19},\\
 u_*&=\frac{(1+3b_*)(1+6b_*-3b_*^2)^{3/2}}
          {12\sqrt{b_*}\,(1+3b_*^2)},\\
 v_*&=\frac{(1-3b_*)(3b_*^2+6b_*-1)^{3/2}}
          {12\sqrt{b_*}\,(1+3b_*^2)}.
 \end{aligned}
\]
\end{theorem}

\begin{proof}
By \eqref{eq:overview-distance-to-ellipse}, it suffices to estimate
the distance from \(C_b\) to \(\partial E_b\) for each \(b\).
First suppose that \(0\le b\le1/4\). The inequality at \(s=3\)
gives \(C_b\subseteq K_b\), where
\[
 K_b=\left\{u+iv\in\mathbb C:
 \frac{3u^2}{(1+2b)^2}+\frac{3v^2}{(1-2b)^2}\le1\right\},
\]
as illustrated in Figure~\ref{fig:coefficient-enclosure}.

\begin{figure}[H]
\centering
\includegraphics[width=\textwidth]{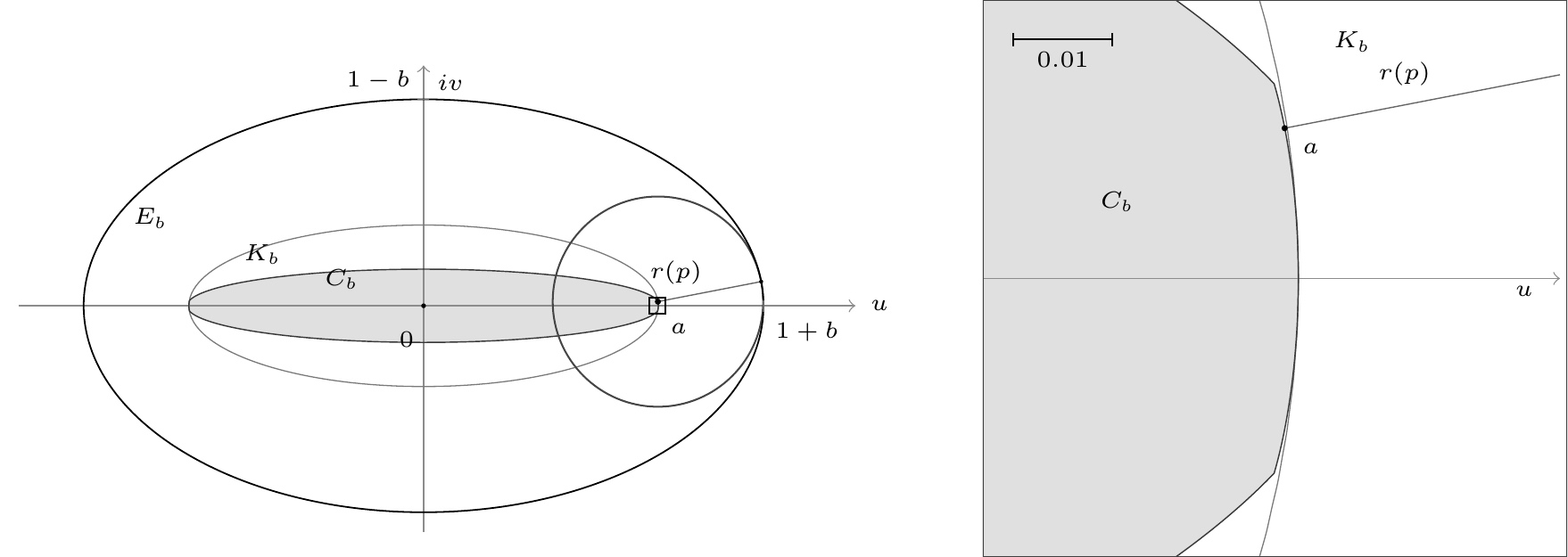}
\caption{The coefficient region \(C_b\) and its containing ellipse
\(K_b\) for \(b=10/41\), together with a magnification of the boxed
region.}
\label{fig:coefficient-enclosure}
\end{figure}

We now have to estimate the distance from \(K_b\) to \(\partial E_b\).
The ellipses \(\partial K_b\) and \(\partial E_b\) are concentric and
have the same axes, but they are not confocal when \(b>0\). Consequently,
a shortest segment joining them need not lie on the \(u\)-axis. In fact,
it lies on this axis precisely when
\(0\leq b\leq(\sqrt3-1)/4\) (see Figure~\ref{fig:auxiliary-ellipse}).

This distance can be studied in two ways. The first is to construct
parallel supporting lines to the two ellipses and minimize the distance
between them over the direction of their common normal. We take the
second route and introduce an intermediate family of auxiliary ellipses
\(Q_\lambda(a)=c\). This method gives the required uniform estimate
directly and will also be useful in the second case, where
\(1/4\le b<1/3\).

Denote the semiaxes of \(E_b\) by
\[
 \alpha=1+b,\qquad \beta=1-b.
\]
Then every point of \(\partial E_b\) can be written as
\[
 w=\alpha x+i\beta y,\qquad x^2+y^2=1.
\]
Fix a point \(a=u+iv\in E_b\). We now use the completion-of-squares
argument underlying the Lagrange multiplier equations for this
constrained nearest-point problem. For every \(\lambda\) different from both
\(\alpha^2\) and \(\beta^2\), and every \(w\in\partial E_b\), the identity
\begin{equation}\label{eq:distance-squares}
\begin{aligned}
 |w-a|^2
 &=\lambda\bigl(1-Q_\lambda(a)\bigr)\\
 &\quad +(\alpha^2-\lambda)
 \left(x-\frac{\alpha u}{\alpha^2-\lambda}\right)^2
 +(\beta^2-\lambda)
 \left(y-\frac{\beta v}{\beta^2-\lambda}\right)^2
\end{aligned}
\end{equation}
holds, where
\begin{equation}\label{eq:Q-lambda}
 Q_\lambda(a)
 =\frac{u^2}{\alpha^2-\lambda}
  +\frac{v^2}{\beta^2-\lambda}.
\end{equation}
The two squares in \eqref{eq:distance-squares} vanish simultaneously
precisely when
\begin{equation}\label{eq:lagrange-system}
 (\alpha^2-\lambda)x=\alpha u,\qquad
 (\beta^2-\lambda)y=\beta v.
\end{equation}
In this case, \(\lambda\) is the corresponding Lagrange multiplier.
If \(0<\lambda<\beta^2\), then both coefficients of the squares in
\eqref{eq:distance-squares} are positive. Hence
\begin{equation}\label{eq:pointwise-distance-lower}
 |w-a|^2\geq\lambda\bigl(1-Q_\lambda(a)\bigr).
\end{equation}
Equality holds precisely when \eqref{eq:lagrange-system} holds.
Thus, for every \(a\in E_b\) and every \(\lambda\) satisfying
\(0<\lambda<\beta^2\),
\begin{equation}\label{eq:distance-to-boundary-lower}
 \operatorname{dist}(a,\partial E_b)^2
 \geq\lambda\bigl(1-Q_\lambda(a)\bigr).
\end{equation}
The identity \eqref{eq:distance-squares}, together with the equality
condition \eqref{eq:lagrange-system}, will be used not only to estimate
this distance from below but also later to identify the critical point.

Our task is now to obtain a uniform lower bound for the distance in
\eqref{eq:distance-to-boundary-lower} as \(a\) ranges over the ellipse
\(K_b\). For this purpose, we choose \(\lambda\) and \(c\) so that
\[
 K_b=\{a\in\mathbb C:Q_\lambda(a)\leq c\}.
\]
In other words, a level ellipse of \(Q_\lambda\) must coincide with
\(K_b\). Matching the semiaxes requires
\[
 c(\alpha^2-\lambda)=\frac{(1+2b)^2}{3},
 \qquad
 c(\beta^2-\lambda)=\frac{(1-2b)^2}{3}.
\]
For \(b>0\), these equations give
\[
 c=\frac23,\qquad \lambda=\frac{1-2b^2}{2},
\]
and the same choice works when \(b=0\). Thus, for the chosen
\(\lambda\), every \(a\in K_b\) satisfies
\(Q_\lambda(a)\leq2/3\). Hence, by
\eqref{eq:distance-to-boundary-lower},
\begin{equation}\label{eq:small-b-bound}
 \operatorname{dist}(a,\partial E_b)^2
 \geq\lambda\left(1-\frac23\right)
 =\frac{1-2b^2}{6}\geq\frac7{48},
 \qquad 0\leq b\leq\frac14.
\end{equation}

\begin{figure}[H]
\centering
\includegraphics[width=.82\textwidth]{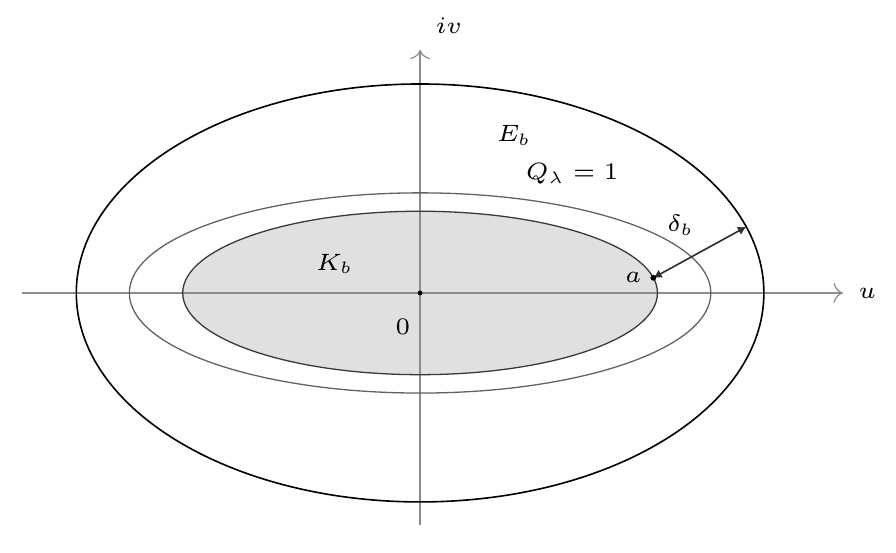}
\caption{The three nested ellipses for \(b=10/41\) and
\(\lambda=(1-2b^2)/2\). The middle ellipse is
\(Q_\lambda=1\), while \(K_b\) is the level region
\(Q_\lambda(a)\leq2/3\). The double arrow from \(a\in\partial K_b\)
has length \(\delta_b=\operatorname{dist}(K_b,\partial E_b)\).}
\label{fig:auxiliary-ellipse}
\end{figure}

Now let \(1/4\le b<1/3\). We seek a sharp bound for this fixed \(b\).
The point where the arc of the ellipse \(s=4\) meets the envelope
of Brannan's family is a natural candidate for the minimum; see
Figure~\ref{fig:coefficient-corner}.

\begin{figure}[H]
\centering
\includegraphics[width=.84\textwidth]{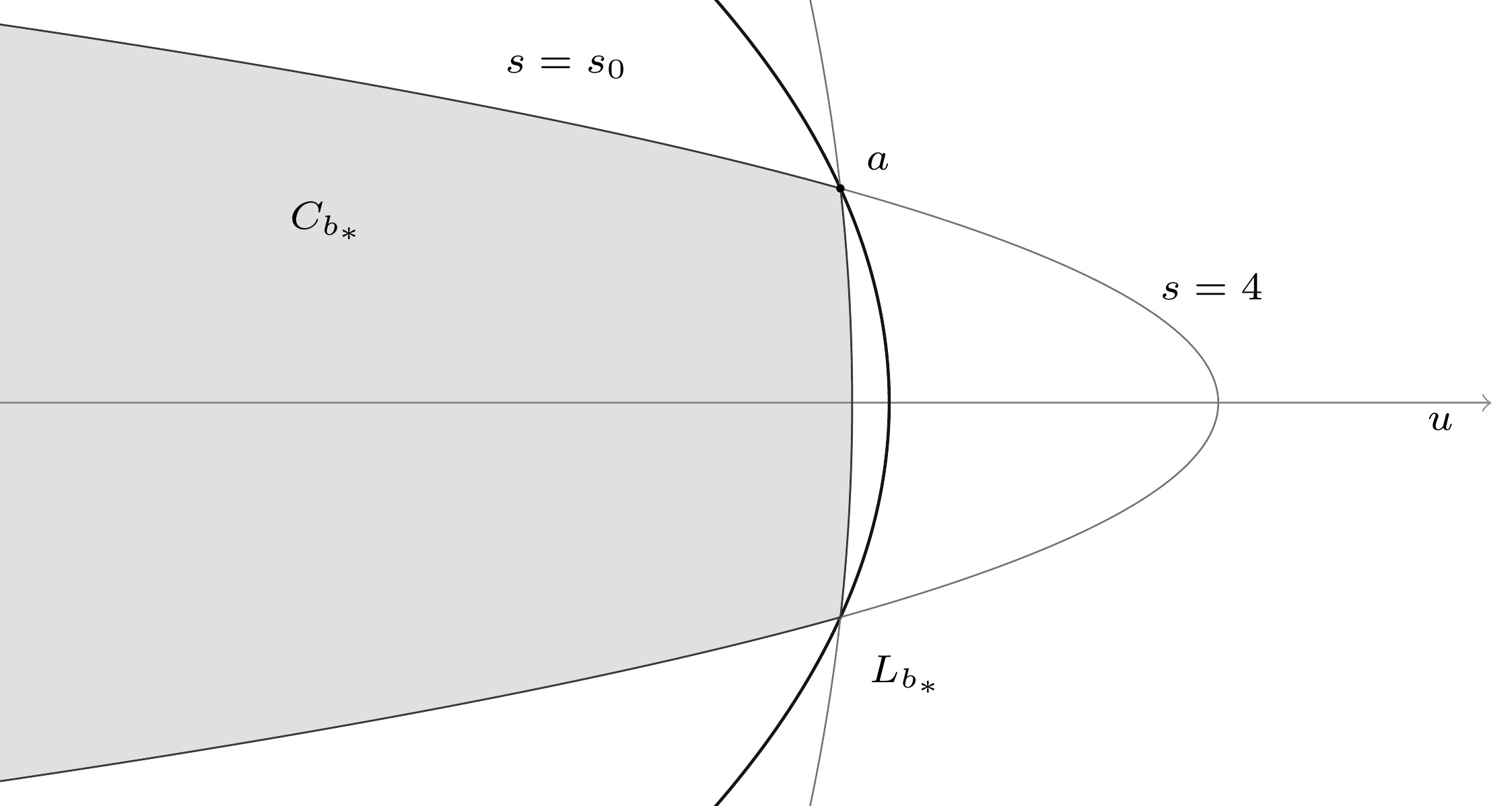}
\caption{The coefficient region \(C_{b_*}\), the two members of
Brannan's family corresponding to \(s=s_0\) and \(s=4\), and the
auxiliary ellipse \(L_{b_*}\) constructed below. Their boundaries
meet at the corner point \(a\).}
\label{fig:coefficient-corner}
\end{figure}

Our next task is to find the two ellipses in Brannan's family whose
boundaries meet at the corner point \(a\). One ellipse is already
known: it corresponds to \(s=4\). The second is unknown; denote its
parameter by \(s_0\). This second ellipse must not merely pass through
\(a\), but must generate the envelope of the family at that point.
Writing \(a=u+iv\), put
\[
 \Phi_b(s;u,v)=
 \frac{s u^2}{(1+b(s-1))^2}
 +\frac{s v^2}{(1-b(s-1))^2}.
\]
Thus \(\Phi_b(s;u,v)=1\) is the boundary of the ellipse with parameter
\(s\). We must therefore find \(u,v,s_0\) satisfying the three
conditions
\[
 \Phi_b(4;u,v)=1,\qquad
 \Phi_b(s_0;u,v)=1,\qquad
 \frac{\partial\Phi_b}{\partial s}(s_0;u,v)=0.
\]

We first solve the two envelope equations
\[
 \Phi_b(s;u,v)=1,\qquad
 \frac{\partial\Phi_b}{\partial s}(s;u,v)=0.
\]
Since
\[
 \frac{\partial}{\partial s}
 \frac{s}{(1\pm b(s-1))^2}
 =\frac{1\mp b(1+s)}{(1\pm b(s-1))^3},
\]
solving for \(u^2\) and \(v^2\) gives
\[
 \begin{aligned}
 u^2(s)&=
 \frac{(1+b(1+s))(1+b(s-1))^3}{4bs^2},\\
 v^2(s)&=
 \frac{(b(1+s)-1)(1-b(s-1))^3}{4bs^2}.
 \end{aligned}
\]
We now impose the remaining condition
\(\Phi_b(4;u(s),v(s))=1\). Direct substitution gives
\[
 \Phi_b(4;u(s),v(s))-1
 =-\frac{(s-4)^2(6b^2s-1-3b^2)(2b^2s-1+b^2)}
 {s^2(1-9b^2)^2}.
\]
The root \(s=4\) gives the ellipse already known, whereas
\[
 \frac{1-b^2}{2b^2}>4
\]
when \(1/4\le b<1/3\). Hence the only parameter of a distinct ellipse
in Brannan's family is
\begin{equation}\label{eq:s-zero}
 s_0=\frac{1+3b^2}{6b^2}.
\end{equation}

Here \(2<s_0\le19/6<4\). Setting \(s=s_0\) in the envelope formulas
gives the first-quadrant intersection of the two ellipses:
\begin{equation}\label{eq:fixed-b-extremizer}
 \begin{aligned}
 u_b&=\frac{(1+3b)(1+6b-3b^2)^{3/2}}
          {12\sqrt b\,(1+3b^2)},\\
 v_b&=\frac{(1-3b)(3b^2+6b-1)^{3/2}}
          {12\sqrt b\,(1+3b^2)},\\
 a_b&=u_b+iv_b,\qquad
 p_b(z)=z+a_bz^2+bz^3.
 \end{aligned}
\end{equation}
Both radicals have positive arguments on the interval under
consideration.

To check that \(a_b\) satisfies every inequality in
\eqref{eq:intro-ellipses}, put \(\eta=b^2\), \(d=b(s-1)\) and
\[
 P(s)=(1-d^2)^2\bigl(1-\Phi_b(s;u_b,v_b)\bigr).
\]
This is a quartic with leading coefficient \(\eta^2\). By construction it
vanishes at \(4\) and has a double zero at \(s_0\). Since
\(P(0)=(1-\eta)^2\),
\[
 P(s)=\eta^2(s-4)(s-s_0)^2(s-s_1),\qquad
 s_1=\frac{(1-\eta)^2}{4\eta^2s_0^2}
    =\frac{9(1-\eta)^2}{(1+3\eta)^2}>4,
\]
where the last inequality follows from \(\eta<1/9\). Thus \(P(s)\ge0\)
on \([0,4]\). Since \(|d|<1\) there, \(a_b\in C_b\) and \(p_b\)
is univalent.

We now determine the distance from \(a_b\) to \(\partial E_b\).
We choose \(0<\lambda<\beta^2\) and \(w_b\in\partial E_b\) so that
both squares in \eqref{eq:distance-squares} vanish. This gives
equality in \eqref{eq:distance-to-boundary-lower} at \(a_b\).
We then show that \(Q_\lambda(a)\le Q_\lambda(a_b)\) for every
\(a\in C_b\), so \(a_b\) minimizes the distance over \(C_b\).

The equality \(\Phi_b(4;u_b,v_b)=1\) shows that
\[
 x_b=\frac{2u_b}{1+3b},\qquad
 y_b=\frac{2v_b}{1-3b}
\]
satisfy \(x_b^2+y_b^2=1\). Thus
\(w_b=\alpha x_b+i\beta y_b\) lies on \(\partial E_b\). We choose
\(\lambda\) so that the Lagrange equations
\eqref{eq:lagrange-system} hold for \(a_b\) and \(w_b\). The first
equation gives
\[
 \lambda=\alpha^2-\frac{\alpha u_b}{x_b}
 =\frac{1-b^2}{2},
\]
and the second gives the same value. Indeed,
\[
 \alpha^2-\lambda=\frac{(1+b)(1+3b)}{2},
 \qquad
 \beta^2-\lambda=\frac{(1-b)(1-3b)}{2}.
\]
Both quantities are positive when \(1/4\le b<1/3\). Hence the two
squares in \eqref{eq:distance-squares} vanish at \(a_b,w_b\).

It remains to prove that \(Q_\lambda\) is largest on \(C_b\) at
\(a_b\). For \(0<s\le4\), denote the left-hand side of Brannan's
inequality by
\[
 \mathcal B_s(a)=A_+(s)u^2+A_-(s)v^2,
 \qquad
 A_\pm(s)=\frac{s}{(1\pm b(s-1))^2}.
\]
Every point of \(C_b\) satisfies \(\mathcal B_s(a)\leq1\), while
\(\mathcal B_{s_0}(a_b)=\mathcal B_4(a_b)=1\). We shall express
\(Q_\lambda\) as a positive linear combination of
\(\mathcal B_{s_0}\) and \(\mathcal B_4\). Such a representation
immediately gives the desired maximum at \(a_b\).

Write \(R(s)=A_+(s)/A_-(s)\). Both coefficients of each form
\(\mathcal B_s\) are positive, so the required combination exists
precisely when
\[
 R(4)<k<R(s_0),\qquad
 k=\frac{\beta^2-\lambda}{\alpha^2-\lambda}
  =\frac{(1-b)(1-3b)}{(1+b)(1+3b)}.
\]
Substituting \(s_0\) from \eqref{eq:s-zero} gives
\[
 \begin{aligned}
 k-R(4)&=\frac{4b(1-3b)}{(1+b)(1+3b)^2},\\
 R(s_0)-k&=\frac{16b(18b^4+15b^2-1)}
                  {(1+b)(1+3b)(1+6b-3b^2)^2}.
 \end{aligned}
\]
The first difference is positive for \(0<b<1/3\). The numerator of
the second equals \(1/32\) at \(b=1/4\) and increases for \(b\ge1/4\).
Thus for \(1/4\le b<1/3\) we can choose positive \(\mu_b,\nu_b\)
such that
\[
 Q_\lambda=\mu_b\mathcal B_{s_0}+\nu_b\mathcal B_4.
\]
The lower threshold for both weights to be positive is
\(b=\sqrt{(\sqrt{33}-5)/12}=0.249092\ldots\), just below \(1/4\).
At this value \(R(s_0)=k\) and the weight of \(\mathcal B_4\) vanishes.

Set \(\gamma_b=Q_\lambda(a_b)\). Since both Brannan inequalities are
equalities at \(a_b\), we have \(\gamma_b=\mu_b+\nu_b\). Hence
\[
 Q_\lambda(a)\leq\gamma_b,\qquad a\in C_b,
\]
and equality forces both Brannan inequalities to be equalities. Thus
the ellipse
\[
 L_b=\{a:Q_\lambda(a)\leq\gamma_b\}
\]
contains \(C_b\) and touches it at \(a_b\), as shown in
Figure~\ref{fig:coefficient-corner}. Combining this inclusion with
\eqref{eq:distance-to-boundary-lower} gives
\[
 \operatorname{dist}(a,\partial E_b)^2
 \geq\lambda(1-\gamma_b),\qquad a\in C_b.
\]
At \(a_b\), the two squares in \eqref{eq:distance-squares} vanish at
\(w_b\), so equality holds. Hence \(p_b\) is extremal for this fixed
value of \(b\).

It remains to compute the distance. Since \(x_b^2+y_b^2=1\),
\[
 \gamma_b=\frac{1+3b}{2(1+b)}\,x_b^2
        +\frac{1-3b}{2(1-b)}\,y_b^2,
 \qquad
 \lambda(1-\gamma_b)=\frac{(1+b)^2}{4}-b x_b^2.
\]
Recalling that \(\eta=b^2\) and substituting \(u_b\) from
\eqref{eq:fixed-b-extremizer}, we obtain the exact squared radius for
this fixed \(b\):
\begin{equation}\label{eq:D}
 r(p_b)^2=D(\eta)
 =\frac{27\eta^3+108\eta^2-9\eta+2}{9(1+3\eta)^2}.
\end{equation}
At \(b=1/3\), apply the bound to \(p(hz)/h\), whose third coefficient
is \(h^2/3\). Uniform convergence on the unit circle as \(h\uparrow1\)
gives \(r(p)^2\ge D(1/9)=4/27\), with equality for \(p_S\).

It remains to minimize \(D(\eta)\) for
\(1/16\leq\eta\leq1/9\). Differentiation gives
\[
 D'(\eta)=
 \frac{27\eta^3+27\eta^2+81\eta-7}{3(1+3\eta)^3}.
\]
The polynomial in the numerator is strictly increasing for
\(\eta\geq0\) and has opposite signs at \(1/16\) and \(1/9\).
Hence it has a unique zero \(\eta_*\) in this interval, determined by
\begin{equation}\label{eq:eta-star}
 27\eta_*^3+27\eta_*^2+81\eta_*-7=0.
\end{equation}
The sign of \(D'(\eta)\) shows that \(D\) decreases before \(\eta_*\)
and increases after it. Thus \(D(\eta_*)\) is the unique minimum.

Put \(q_*=D(\eta_*)\). Reducing this expression by
\eqref{eq:eta-star} gives
\[
 q_*=\frac{1+15\eta_*}{15+9\eta_*},
 \qquad
 \eta_*=\frac{15q_*-1}{15-9q_*}.
\]
Substitution in \eqref{eq:eta-star}, followed by clearing the
denominator, gives
\[
 G(q_*)=0,
 \qquad
 G(q)=27q^3-54q^2+63q-8.
\]
Moreover,
\[
 G'(q)=81(q-2/3)^2+27>0,
\]
so \(q_*\) is the unique real zero of \(G\). Consequently,
\[
 r_*=\sqrt{q_*}=\sqrt{D(\eta_*)}=0.378591\ldots,
 \qquad
 b_*=\sqrt{\eta_*}=0.2896\ldots.
\]
Thus \(b_*>1/4\). For \(b\le1/4\), the preceding lower bound
\(\sqrt{7/48}\) is strictly larger than \(r_*\), since
\[
 G(7/48)=\frac{503}{4096}>0.
\]
The inequality is also strict at \(b=1/3\), since
\(1/9\ne\eta_*\).
Writing \(q=r_*^2=q_*\), the identity
\[
 (15-9q)\left(q^2-\frac q3+\frac19\right)-(15q-1)
 =-\frac13G(q)
\]
then gives
\[
 b_*=\sqrt{r_*^4-\frac{r_*^2}{3}+\frac19}.
\]
Putting \(t=3r_*^2-2\) reduces \(G(r_*^2)=0\) to
\[
 t^3+9t+18=0,
\]
which yields \eqref{eq:r-star}.

Finally, if \(r(p)=r_*\), the strict bounds at \(b\le1/4\) and \(b=1/3\)
and the unique minimum of \(D\) force \(b=b_*\).
At a nearest point of \(\partial E_b\), both squares in
\eqref{eq:distance-squares} must vanish and
\(Q_\lambda(a)=\gamma_b\).
Hence
\[
 \Phi_b(s_0;u,v)=\Phi_b(4;u,v)=1.
\]
These two equations are linear in \(u^2,v^2\). Their coefficient
ratios \(R(s_0)\) and \(R(4)\) are distinct, so their only solution is
\(u^2=u_{b_*}^2,\ v^2=v_{b_*}^2\). The normalized rotation
\(p(z)\mapsto-p(-z)\) changes \(a\) to \(-a\); conjugation changes
\(a\) to \(\bar a\). Together they give all four choices of signs.

\end{proof}

To locate the cusp of \(p_*(\partial\Delta)\), take \(b=b_*\) and
put \(\zeta=-x_b+i y_b\). Then \(|\zeta|=1\), and the formulas for
\(x_b,y_b\) give \(2a_b+\bar\zeta+3b\zeta=0\). It follows that
\(p_*'(\zeta)=0\). The product of the two zeros of \(p_*'\) is
\(1/(3b)\), so the other critical point has modulus \(1/(3b)>1\).
Moreover, \(p_*(\zeta)/\zeta^2=a_b-w_b\), hence
\(|p_*(\zeta)|=|a_b-w_b|=r_*\).

Write \(\gamma(t)=p_*(\zeta e^{it})\) and \(h=3b\zeta^2\). Then
\(\gamma'(0)=0\), \(\gamma''(0)=\zeta(1-h)\), and
\(\gamma'''(0)=i\zeta(3-5h)\). Since \(\Re h=3b(2x_b^2-1)\),
\eqref{eq:fixed-b-extremizer} gives
\[
 \Im\bigl(\overline{\gamma''(0)}\,\gamma'''(0)\bigr)
 =3-8\Re h+5|h|^2
 =\frac{(21b^2-1)^2(3b^2+5)}{3(1+3b^2)^2}>0.
\]
Thus the quadratic and cubic terms in the Taylor expansion point in
linearly independent directions, and \(\gamma\) has an ordinary cusp at
\(p_*(\zeta)\).

\section*{Acknowledgement}

The author is grateful to Dmitry Dmitrishin for drawing his attention
to this problem. The author has also freely used results of
mathematicians published in open sources, as well as various large
language models.


\begin{thebibliography}{99}

\bibitem{Alexander}
J. W. Alexander,
\href{https://www.jstor.org/stable/2007212}
{\emph{Functions which map the interior of the unit circle upon simple
regions}},
Ann. of Math. (2) \textbf{17} (1915), no.~1, 12--22.

\bibitem{Bieberbach}
L. Bieberbach,
\href{https://archive.org/details/sitzungsberichte1916deut/page/940/mode/2up}
{\emph{\"Uber die Koeffizienten derjenigen Potenzreihen, welche eine
schlichte Abbildung des Einheitskreises vermitteln}},
Sitzungsber. Preuss. Akad. Wiss. Berlin (1916), 940--955.

\bibitem{Brannan}
D. A. Brannan,
\href{https://doi.org/10.1112/S0025579300003764}
{\emph{Coefficient regions for univalent polynomials of small degree}},
Mathematika \textbf{14} (1967), 165--169.

\bibitem{CowlingRoyster}
V. F. Cowling and W. C. Royster,
\href{https://doi.org/10.2307/2035305}
{\emph{Domains of variability for univalent polynomials}},
Proc. Amer. Math. Soc. \textbf{19} (1968), 767--772.

\bibitem{deBranges}
L. de Branges,
\href{https://doi.org/10.1007/BF02392821}
{\emph{A proof of the Bieberbach conjecture}},
Acta Math. \textbf{154} (1985), 137--152.

\bibitem{Dimitrov}
D. K. Dimitrov,
\href{https://www.dcce.ibilce.unesp.br/~dimitrov/papers/main.pdf}
{\emph{Extremal positive trigonometric polynomials}},
in \emph{Approximation Theory: A Volume Dedicated to Blagovest Sendov},
B. Bojanov, ed., DARBA, Sofia, 2002, 136--157.

\bibitem{DDS}
D. Dmitrishin, K. Dyakonov and A. Stokolos,
\href{https://arxiv.org/abs/1812.08311}
{\emph{Univalent polynomials and Koebe's one-quarter theorem}},
Anal. Math. Phys. \textbf{9} (2019), 991--1004.

\bibitem{Duren}
P. L. Duren,
\href{https://link.springer.com/book/9780387907956}
{\emph{Univalent Functions}},
Grundlehren der mathematischen Wissenschaften \textbf{259},
Springer-Verlag, New York, 1983.

\bibitem{Koebe}
P. Koebe,
\href{http://resolver.sub.uni-goettingen.de/purl?GDZPPN002501198}
{\emph{\"Uber die Uniformisierung beliebiger analytischer Kurven}},
Nachr. Kgl. Ges. Wiss. G\"ottingen, Math.-Phys. Kl. (1907), 191--210.

\bibitem{RahmanWaniurski}
Q. I. Rahman and J. Waniurski,
\href{https://doi.org/10.4153/CJM-1980-001-0}
{\emph{Coefficient regions for univalent trinomials}},
Canad. J. Math. \textbf{32} (1980), 1--20.

\bibitem{Suffridge}
T. J. Suffridge,
\href{https://doi.org/10.1112/jlms/s1-44.1.496}
{\emph{On univalent polynomials}},
J. London Math. Soc. \textbf{44} (1969), 496--504.

\end{thebibliography}
\end{document}